\documentclass[number,citesort,seceqn,dvips]{arxbj}
\usepackage{mathbh}
\usepackage[vtexbibl]{}

\aid{0}
\volume{16}
\issue{3}
\pubyear{2010}
\firstpage{605}
\lastpage{613}
\doi{10.3150/09-BEJ225}

\makeatletter

\newtheorem{Theorem}{Theorem}[section]
\newtheorem{Lemma}[Theorem]{Lemma}

\def\eqref#1{(\ref{#1})}

\makeatother

\begin{document}
\begin{frontmatter}

\title{Sharper lower bounds on the performance of the empirical risk
minimization algorithm}
\runtitle{Lower bounds for ERM}

\begin{aug}
\author[1]{\fnms{Guillaume} \snm{Lecu\'{e}}\corref{}\thanksref{1}\ead[label=e1]{lecue@latp.univ-mrs.fr}}
\and
\author[2]{\fnms{Shahar} \snm{Mendelson}\thanksref{2}\ead[label=e2]{shahar.mendelson@anu.edu.au}}
\runauthor{G. Lecu\'{e} and S. Mendelson}
\address[1]{CNRS, LATP, Marseille 13000, France. \printead{e1}}
\address[2]{Centre for Mathematics and Its Applications, The
Australian National University, Canberra, ACT 0200, Australia and
Department of Mathematics, Technion, I.I.T., Haifa 32000, Israel.\\
\printead{e2}}
\end{aug}

\received{\smonth{10} \syear{2008}}
\revised{\smonth{5} \syear{2009}}

%
\begin{abstract}
We present an argument based on the multidimensional and the
uniform central limit theorems, proving that, under some geometrical
assumptions between the target function $T$ and the learning class $F$, the
excess risk of the empirical risk minimization algorithm is lower
bounded by
\[
\frac{\mathbb{E}\sup_{q \in Q}G_{q}}{\sqrt{n}} \delta,
\]
where $(G_{q})_{q \in Q}$ is a canonical Gaussian process
associated with $Q$ (a well chosen subset of $F$) and $\delta$ is a parameter
governing the oscillations of the empirical excess risk function over a
small ball in $F$.
\end{abstract}

%
\begin{keyword}
\kwd{empirical risk minimization}
\kwd{learning theory}
\kwd{lower bound}
\kwd{multidimensional central limit theorem}
\kwd{uniform central limit theorem}
\end{keyword}

\end{frontmatter}

\section{Introduction}
In this note, we study lower bounds on the empirical minimization
algorithm. To explain the basic setup of this algorithm, let
$(\Omega,\mu)$ be a probability space and set $X$ to be a random
variable taking values in $\Omega$, distributed according to $\mu$.
We are interested in the {\it function learning} (noiseless)
problem, in which one observes $n$ independent random variables
$X_1,\ldots,X_n$, distributed according to $\mu$, and the values
$T(X_1),\ldots,T(X_n)$ of an unknown target function $T$.

The goal is to construct a procedure that uses the data
$D=(X_i,T(X_i))_{1\leq i \leq n}$ with a {\it risk} as close as
possible to the best one in $F$. That is, we want to construct a
statistic $\hat{f}_n$ such that for every $n$, with high
$\mu^n$-probability,
%
\begin{equation}\label{Equa:InegaliteOracle}
R(\hat{f}|D)\leq\inf_{f\in F} R(f) + r_n(F),
\end{equation}
where the risk of $f$ is defined by $R(f)=\mathbb{E}\ell(f(X),T(X))$ and
$\ell\dvtx\mathbb{R}^2 \to\mathbb{R}$ is the loss function that
measures the pointwise
error between $T$ and $f$. The residue $r_n(F)$ somehow captures the
complexity or richness of the class $F$ and the risk of a statistic
$\hat{f}$ is the conditional\vspace{1pt} expectation $R(\hat{f}|D)=\mathbb{E}
(\ell(\hat{f}(X),T(X))|D)$.

It is well known (see, e.g., \cite{VapBook98}) that if the class $F$
is not too large, for example, if it satisfies some kind of uniform central
limit theorem, $T$ is bounded by $1$ and $\ell$ is reasonable, then there
are upper bounds on $r_n(F)$ that are of the form $\sqrt{{\operatorname{Comp}(F)}/{n}}$,
where $\operatorname{Comp}(F)$ is a complexity term that is independent of
$n$. The algorithm that is used to produce the function $\hat{f}$ is
the empirical risk minimization algorithm, in which one chooses a
function in $F$ that minimizes the empirical risk function
$f\longmapsto\sum_{i=1}^n \ell(f,T)(X_i)$ in $F$.

There is a well developed theory concerning ways in which the complexity
term may be controlled, using various parameters associated with the
geometry of the class (cf.~\cite{vdVWBook96,VapBook98,DudleyBook99,TalBook05}
and references therein).
It turns out that this type of error
rate, $\sim1/\sqrt{n}$, is very pessimistic in many cases. In fact,
if the class is small enough, then, under some structural
assumptions (see, e.g., \cite{BM06}), $r_n(F)$ can be much
smaller -- of the order of ${\operatorname{Comp}(F)}/{n}$.

In this note, we are going to focus on ``small classes'' in which
empirical minimization performs poorly, despite the size of the
class. Recently, it has been shown (cf.~\cite{M08}) that under mild
assumptions on
$\ell$ and $F$, if there is more than a single function in
\[
V := \Bigl\{\ell(f,T) \dvt \mathbb{E}\ell(f,T) = \inf_{f \in F} \mathbb
{E}\ell(f,T) \Bigr\},
\]
then the
following holds: for every $n$ large enough, there will be a
perturbation $T_n$ of $T$ (with respect to the $L_\infty$-norm) for
which $\mathbb{E}\ell(\cdot, T_n)$ has a unique minimizer in $F$,
but where the
empirical minimization algorithm performs poorly trying to predict
$T_n$ on samples of cardinality $n$. To be more exact, relative to
the target $T_n$, with $\mu^n$-probability at least $1/12$,
%
\begin{equation} \label{eq:01}
R(\hat{f}|D) \geq\inf_{f\in F} R(f) + \frac{c}{\sqrt{n}},
\end{equation}
where $c$ is a constant depending only on $F$.

Although it is reasonable to expect that the larger the set $V$ is, the
more likely it
is that the empirical minimization algorithm will perform poorly, this
does not follow from the analysis in \cite{M08}. Therefore, our goal
here is to provide a bound on the constant $c$ in \eqref{eq:01} that
does take into account the complexity of the set of minimizers
$V$.

Just as in \cite{M08}, our method of analysis can be applied to a
wide variety of losses. However, for the sake of simplicity, we will
only present here what is arguably the most important case -- that in
which the risk is measured relative to the squared loss,
$\ell(x,y)=(x-y)^2$.

To explain our result, we need several definitions from empirical
processes theory. Other standard notions we require from the theory
of Gaussian processes can be found in \cite{DudleyBook99}.

For every set $F \subset L_2(\Omega,\mu)$, let $\{G_f \dvt f \in F\}$ be
the canonical Gaussian process indexed by $F$ (i.e., with the
covariance structure $\mathbb{E}G_t G_s = \langle s,t \rangle$) and
set $H(F)=\mathbb{E}
\sup_{f \in F} G_f$ -- the expectation of the supremum of the
Gaussian process indexed by $F$. Also, for every integer $n$ and
$\delta$, let
\[
\operatorname{osc}_n(F,\delta) := \frac{1}{\sqrt{n}} \mathbb{E}\sup_{\{f,h
\in
F : \|f-h\| \leq\delta\}} \Biggl|\sum_{i=1}^n g_i(f-h)(X_i) \Biggr|,
\]
where $(g_i)_{i=1}^n$ are standard, independent Gaussian random
variables and $(X_i)_{i=1}^n$ are independent, distributed according
to $\mu$. It is well known that if $F$ is a class consisting of
uniformly bounded functions, then it is a $\mu$-Donsker class if and
only if for every $\delta>0$, $\operatorname{osc}_n(F,\delta)$ tends to $0$
as $n$ tends to infinity (cf.~\cite{DudleyBook99}, page 301). For any
$f \in F$, let
\[
\operatorname{osc}_n(F,f,\delta) := \frac{1}{\sqrt{n}} \mathbb{E}\sup_{\{h
\in
F : \|f-h\| \leq\delta\}} \Biggl|\sum_{i=1}^n g_i(f-h)(X_i) \Biggr|,
\]
that is, the oscillation in a ball around $f$. The quantity ${\rm
osc}_n(F,f^*,\delta)$
is a natural upper bound for some intrinsic quantity of the problem we
study here (cf.~Lemma~\ref{lemma:sym}).

Let $V$ be as above -- the set of loss functions $\ell(f,T)$ that
minimize the risk in $F$ -- select $f^* \in F$ for which $\ell(f^*,T)
\in V$ and consider the following subset of excess loss functions:
\[
Q:=\{\ell(f,T)-\ell(f^*,T) \dvt \ell(f,T) \in V\}.
\]

It turns out that the desired constant in (\ref{eq:01}) can be bounded
from below by two
parameters: the expectation of the supremum of the canonical
Gaussian process indexed by $Q$ and the oscillation around $f^*$. In
particular, if $Q$ is a rich set and one of the minimizers of $f \to
\mathbb{E}\ell(f,T)$ is isolated, then for any $n$ large enough, the error
of the empirical minimizer with respect to a wisely selected target
(denoted by $T_{\lambda_n}$ in what follows) which is a perturbation
of $T$ will be at least $\sim H(Q)/\sqrt{n}$. The core idea of this
work is that a small, wisely chosen perturbation of a target function
$T$ with multiple oracles (functions achieving $\min_{f\in F}\mathbb
{E}\ell
(t,T)$) is badly estimated by the empirical risk minimization procedure
(for further discussion of this fact, we refer the reader to \cite{M08}).

Although the general philosophy of the proof presented here is
similar to the proof from \cite{M08}, it is much simpler. And, in
fact, it seems that the method used in the proof from \cite{M08}
cannot be directly extended to obtain the sharper estimate on the
constant as we do here. Naturally, this result recovers the previous
estimates on lower bounds for the empirical risk minimization
algorithm from \cite{LBW98,LecSub07,MN06,K06}

Next, a word about notation. Throughout, all absolute constants will
be denoted by $c, c_1$ and $C, C_1, $ etcetera. Their values may change
from line to line.

If $\mathbb{E}\ell(\cdot,T)$ has a unique minimizer in $F$, then we denote
it by
$f^*$. If the minimizer is not unique, then we will fix one function in
the set of minimizers and denote it by $f^*$. For every $f \in F$,
let $\mathcal{L}(f)=\ell(f,T)-\ell(f^*,T)$ be the excess loss function
associated with the target $T$. For every $0 < \lambda\leq1$, set
$T_\lambda=(1-\lambda) T + \lambda f^*$ and define
$\mathcal{L}_\lambda(f)=\ell(f,T_\lambda)-\ell(f^*,T_\lambda)$.
It is
standard to verify (cf.~\cite{M08} or Theorem~\ref
{cor:OtherMinimizerInCorona} in what follows) that $f^*$ is a minimizer
of $\mathbb{E}\ell(\cdot,
T_\lambda)$ and that under mild convexity assumptions on $\ell$
that clearly hold if $\ell$ is the squared loss, it is the unique
minimizer in $F$ of $f \to\mathbb{E}\ell(f,T_\lambda)$.

If $X_1,\ldots,X_n$ is an independent sample selected according to
$\mu$, we set $P_n f = n^{-1}\times\break \sum_{i=1}^n f(X_i)$ and let $P f =
\mathbb{E}
f$. Thus, $\mathbb{E}\sup_{f \in F} |(P_n-P)(f)|$ is the expectation
of the
supremum of the empirical process indexed by $F$. Finally, when the
target function is $T_\lambda$, we will denote the function
produced by the empirical risk minimization algorithm by
$\hat{f}_{\lambda}$ -- which is one element of the set
$\operatorname{Arg}\min_{f\in F} P_n \ell(f,T_\lambda).$

Finally, if $E$ is a normed space, we denote its unit ball by $B(E)$,
the inner product of $L_2(\mu)$ will be denoted by $\langle\cdot
,\cdot \rangle$ and the corresponding norm by $\|\cdot\|$.

Let us now formulate our main result.

\begin{Theorem} \label{thm:main}
Let $F \subset L_2(\mu) \cap B(L_\infty)$, which is $\mu
$-pre-Gaussian (cf.~\cite{DudleyBook99}),
and assume that $T \in B(L_\infty)$. Set $\ell$ to be the squared
loss and put $Q=\{\mathcal{L}(f) \dvt f \in F,  \mathbb{E}\mathcal{L}(f)=0\}$.

There exist some absolute constants $C_1$ and
$C_2$ and an integer $N(F)$ for which the following holds. For
every $n \geq N(F)$, with $\mu^n$-probability at least $C_1$,
\[
\mathbb{E}\mathcal{L}_{\lambda_n}(\hat{f}_{\lambda_n}) \geq C_2 \frac
{H(Q)}{\sqrt{n}} \delta^2 \|T-f^*\|,
\]
where $\delta$ is such that for every integer $n \geq N(F)$, ${\rm
osc}_n(F,f^*,\delta) \leq C_2 H(Q)/\sqrt{n}$ and $\lambda_n=C_2
H(Q)/\sqrt{n}$.
\end{Theorem}

Thus, two parameters control the behavior of the constant in (\ref{eq:01}):
the complexity of the set of excess loss functions of the
oracles of $T$ and the parameter $\delta$. When one of the oracles
$f^*$ of $T$ is isolated, one can take $\delta$ as an absolute
constant. This leads to a lower bound of the order of $H(Q)/\sqrt{n}$,
which is optimal in the sense that an upper bound can be obtained of
the order of $H(Q_0)/\sqrt{n}$ for some set $Q_0$ such that $Q\subset
Q_0 \subset\mathcal{L}_F$ (see, e.g., \cite{BM06} or \cite{K06}). In
other settings, the lower bound obtained in Theorem \ref{thm:main} may
fail to match exactly with an upper bound. For instance, in settings
where the oscillation function $\operatorname{osc}_n(F,f^*,\cdot)$ of all the
oracles $f^*$ of $T$ decreases to zero very slowly and at the same
convergence rate, the factor $\delta^2$ should break down the lower
bound, whereas it seems that it should not appear in the lower bound.
From a technical point of view, this comes from the fact that we did
not take into account the complexity ``around'' the points in $Q^\prime
$ (cf.~Theorem~\ref{thm:gauss-proc} and equation~(\ref{eq:Final1}) in
what follows).

Finally, the noiseless model considered here is the worst case scenario
to prove the lower bound. Indeed, adding some noise to the target
function would increase the lower bound.

\section{The lower bound}

The core of the proof is to find a set that can ``compete'' with a
set $B_r=\{f\in F\dvt \mathbb{E}\mathcal{L}_\lambda(f)\leq r\}$ that
contains $f^*$, in
the sense that the empirical excess risk function
\[
\mathcal{E}_n\dvtx f\in F\longmapsto\frac{1}{n}\sum_{i=1}^n\mathcal
{L}_\lambda(f)(X_i)
\]
will be more negative on the set than it can possibly be on
$B_r$. Once this task is achieved, it is obvious that the empirical
risk minimization algorithm will produce a function $\hat{f}_\lambda$ which
is outside $B_r$ and, thus, with a certain probability,
\begin{eqnarray*}
\mathbb{E}[ \mathcal{L}_{\lambda}(\hat{f}_{\lambda})|D]> r.
\end{eqnarray*}

Hence, the proof consists of two parts. First, we will show that the
empirical excess risk function $\mathcal{E}_n$ is likely to be very negative
on $Q$ and we will then find some $r$ on which the oscillations in
$B_r$ are small.

The first result we need is the following lower estimate on the
expectation of the excess loss relative to the target
$T_\lambda=(1-\lambda) T +\lambda f^*$, according to the distance of
$f$ from $f^*$. This proposition is based on the fact that the
functional $(f,g)\longmapsto\mathbb{E}\ell(f,g)$ inherits a strong convex
structure from the norm and was proven in \cite{M08} in a far more
general situation.

\begin{Theorem}\label{cor:OtherMinimizerInCorona}
Let $D=\sup_{f\in F}\|T-f\|$ and $\rho=\|T-f^*\|$. There exists an
absolute constant $c$ such that for any function $f\in F$,
if $0\leq\lambda\leq1/2,$ $r>0$ and
\begin{eqnarray*}
\frac{r}{\lambda} \leq c \frac{\rho}{D} \|f-f^*\|^2,
\end{eqnarray*}
then
\[
r\leq\mathbb{E}\mathcal{L}_\lambda(f).
\]
\end{Theorem}

Recall that $Q = \{ \mathcal{L}(f) \dvt \mathbb{E}\mathcal{L}(f) = 0,
f\in F \}$ is the set of
excess loss functions associated with the true minimizers of $f
\to\mathbb{E}\ell(f,T)$ in $F$. We will show that if $Q^\prime
\subset Q$ is a
finite set, then for $n$ large enough, with a non-trivial
$\mu^n$-probability there will be some $\mathcal{L}(f) \in Q^\prime$ for
which the empirical error $P_n \mathcal{L}_{\lambda_n}(f)$ is very negative
(for a well chosen $\lambda_n$).

\begin{Theorem} \label{thm:gauss-proc}
There exist constants $c_1, c_2$ and $c_3$, depending only on the
$L_\infty(\mu)$-diameter of $F \cup\{T\}$, for which the following
holds. If $Q^\prime$ is a finite subset of $Q$ that contains $0$,
then there exists an integer $n_0=n_0(Q^\prime)$ such that for every
integer $n \geq n_0$, with $\mu^n$-probability at least $c_1$,
\[
\inf_{\mathcal{L}(f) \in Q^\prime} \frac{1}{n} \sum_{i=1}^n
(\mathcal{L}_{\lambda_n}(f))(X_i) \leq-c_2\frac{H(Q^\prime)}{\sqrt{n}},
\]
where $\lambda_n=c_3{H(Q^\prime)}/{\sqrt{n}}$ and $H(Q^\prime
)=\mathbb{E}
\sup_{q \in Q^\prime} G_q$ is the expectation of the canonical
Gaussian process associated with $Q^\prime$.
\end{Theorem}

\begin{pf} Let $M=|Q^\prime|$ and recall that each $q \in
Q^\prime=\{q_1,\ldots,q_M\}$ has mean zero. Consider the random
vector $U=(q_1(X),\ldots,q_M(X)) \in\mathbb{R}^M$ and let
$(U_i)_{i=1}^\infty$
be independent copies of $U$ (i.e.,
$U_i=(q_1(X_i),\ldots,q_M(X_i))$). By the vector-valued central
limit theorem (see,~e.g., \cite{DudleyBook99}), $n^{-1/2}
\sum_{i=1}^n U_i$ converges weakly to the canonical Gaussian process
indexed by $Q^\prime$, which we denote by $G$. Fix $t \leq0$ and
$0<c<1$, to be given later, for which
\[
A_t = \{x \in\mathbb{R}^M \dvt \forall 1 \leq j \leq M,    x_j > t\}
\]
is such that $p:=\operatorname{Pr}(G \in A_t)\leq c$. Set $n_0=n_0(t,c)$ to be
such that for $n \geq n_0$,
\[
\Biggl| \operatorname{Pr} (G \in A_t) - \operatorname{Pr} \Biggl(n^{-1/2}\sum_{i=1}^n U_i \in
A_t\Biggr) \Biggr| \leq\frac{1-p}{2},
\]
which clearly exists by weak convergence. Since
\begin{eqnarray*}
\operatorname{Pr} \Biggl(\exists 1 \leq j \leq M  \dvt  n^{-1/2} \sum_{i=1}^n\langle
U_i,e_j \rangle \leq t\Biggr) &=& 1 - \operatorname{Pr} \Biggl( n^{-1/2} \sum_{i=1}^n
U_i \in
A_t\Biggr)
\\
&\geq&\frac{1-p}{2}\geq\frac{1-c}{2}=:c_1>0,
\end{eqnarray*}
it follows that, with probability at least $c_1$,
\[
\inf_{q \in Q^\prime} \frac{1}{n} \sum_{i=1}^n q(X_i) \leq
\frac{t}{\sqrt{n}}.
\]
It remains to show that one may take $t=-(\mathbb{E}\sup_{q \in
Q^\prime}
G_q)/4$. Indeed, by the symmetry of the Gaussian process, it follows
that (for this choice of $t$)
\begin{eqnarray*}
p=\operatorname{Pr}(G\in A_t)=\operatorname{Pr} \Bigl(\sup_{q\in Q^{\prime}} G_q < \Bigl(\mathbb{E}\sup
_{q\in
Q^{\prime}} G_q\Bigr)\big/4\Bigr).
\end{eqnarray*}
Let $Z=\sup_{q\in Q^\prime}G_q$ and $\sigma^2=\sup_{q\in
Q^{\prime}}\mathbb{E}G_q^2$. Since $0 \in Q^\prime$, it follows that
if $\mathbb{E}
Z=0$, then it is
clear that $p=1/2$. Otherwise, using the concentration property of
$Z$ around its mean (see, e.g., \cite{vdVWBook96}) and since $\sigma
\leq c_0 \mathbb{E}Z$ (where $c_0$ is an absolute constant), there exists
an absolute constant $A>0$ such that
\begin{eqnarray*}
\mathbb{E}\bigl[Z\mathbh{1}_{[Z\geq\mathbb{E}Z+A\sigma]}
\bigr]\leq(\mathbb{E}Z)/4.
\end{eqnarray*}
Therefore,
\begin{eqnarray*}
\mathbb{E}Z&=&\mathbb{E}\bigl(Z\bigl(\mathbh{1}_{[Z\leq(\mathbb
{E}Z)/4]}+\mathbh{1}_{[(\mathbb{E}Z)/4\leq Z \leq\mathbb{E}
Z +A\sigma]}+\mathbh{1}_{[Z\geq\mathbb{E}Z +A\sigma]} \bigr)
\bigr)\\
&\leq&(\mathbb{E}Z)/2+(\mathbb{E}Z)(1+c_0A)\operatorname{Pr}\bigl((\mathbb
{E}Z)/4\leq Z\bigr).
\end{eqnarray*}
Thus, $\operatorname{Pr}((\mathbb{E}Z)/4\leq Z)\geq[2(1+c_0A)]^{-1}$ and so
$p\leq1-[2(1+c_0A)]^{-1}:=c$ (which is an absolute constant), implying
that, with probability greater
than $c_1$,
\[
\inf_{\mathcal{L}(f) \in Q^\prime} \frac{1}{n} \sum_{i=1}^n
(\mathcal{L}(f))(X_i)
\leq-c_2\frac{\mathbb{E}\sup_{q \in Q^\prime} G_q}{\sqrt{n}}.
\]

Next, observe that for small values of $\lambda$ (as we will have in
our construction), $\mathcal{L}(f)$ is a good approximation of
$\mathcal{L}_\lambda(f)$ with respect to the $L_\infty(\mu)$-norm. Indeed,
$\mathcal{L}_\lambda(f)=\ell(f,T_\lambda)-\ell(f^*,T_\lambda)$ and
$\mathcal{L}(f)=\ell(f,T)-\ell(f^*,T)$; hence, for every $f \in F$,
\begin{eqnarray*}
\|\mathcal{L}_\lambda(f)-\mathcal{L}(f)\|_\infty &\leq&
\|\ell(f,T_\lambda) - \ell(f,T)\|_\infty+ \|\ell(f^*,T_\lambda) -
\ell(f^*,T)\|_\infty
\\
&\leq& 2\|\ell\|_{\rm lip}\|T-T_\lambda\|_\infty=2\lambda
\|\ell\|_{\rm lip} \|T-f^*\|_\infty\leq c_3\lambda.
\end{eqnarray*}
Thus, if one selects $\lambda_n = (c_2/(2c_3)) n^{-1/2} \mathbb
{E}\sup_{q
\in Q^\prime} G_q$, then, with probability greater than $c_1$,
\[
\inf_{\mathcal{L}(f) \in Q^\prime} P_n \mathcal{L}_{\lambda_n}(f)
\leq
-c_2\frac{\mathbb{E}\sup_{q \in Q^\prime} G_q}{2\sqrt{n}}.
\]
\upqed\end{pf}

Fix a finite set $Q^\prime\subset Q$ for which $H(Q^\prime) \geq
H(Q)/2$ and $0\in Q^\prime$. Clearly, such a set exists because $Q$ is
a pre-Gaussian as
a subset of the pre-Gaussian class $\{\mathcal{L}(f) \dvt f \in F\}$. Let
$V^\prime=\{f\in F \dvt \mathcal{L}(f) \in Q^\prime\}$.


Recall that a bounded class of functions $F$ is $\mu$-\textit
{Donsker} if and
only if for every $u>0$, there exist $\delta>0$ and an integer $n_0$
such that for every $n \geq n_0, \operatorname{osc}_n(F,\delta)\leq u$.
Also,\vadjust{\goodbreak} note that $\operatorname{osc}_n(F,f^*,\delta) \leq{\rm
osc}_n(F,\delta)$. Let $u=\eta H(Q^\prime)$, where $\eta$ is an
absolute constant, to be fixed later, and set $\delta$ and $n_1$ to be
such that for $n \geq n_1$,
%
\begin{equation} \label{eq:delta}
\operatorname{osc}_n(F,f^*,\delta) \leq\eta H(Q^\prime)
\end{equation}
(such $\delta$ and $n_1$ necessarily exist because $F$ is $\mu$-Donsker).

The next lemma is standard and follows from a symmetrization
argument combined with Slepian's lemma. Its proof may be found in, for
example, \cite{M08}.
\begin{Lemma} \label{lemma:sym}
There exists an absolute constant $c$ for which the following holds.
For any $F^\prime\subset F$ such that $f^*\in F^\prime$ and any
$0\leq\lambda\leq1$,
\[
\mathbb{E}\sup_{f \in F^\prime} |(P-P_n)(\mathcal{L}_\lambda
(f))| \leq c
\mathbb{E}\sup_{f \in F^\prime} \Biggl|\frac{1}{n} \sum_{i=1}^n
g_i(f-f^*)(X_i)\Biggr|,
\]
where $(g_i)_{i=1}^n$ are independent, standard Gaussian variables.
\end{Lemma}

We are now ready to control the oscillation of the empirical excess
risk function in the set $B_r=\{f\in F\dvt \mathbb{E}\mathcal{L}_\lambda
\leq r\}$.
\begin{Theorem} \label{thm:UCLT}
Let $c_1$, $c_2$ and $\lambda_n$ be defined as in
Theorem~\ref{thm:gauss-proc}, and let $\delta$ and $n_1$ be as
above. There exists an absolute constant $c_3$ such that for any
integer $n\geq n_1$, with $\mu^n$-probability at least $1-c_1/2$,
\begin{eqnarray*}
\inf_{\{f\in F: \mathbb{E}\mathcal{L}_{\lambda_n}(f) \leq r_{n}\}} P_n
\mathcal{L}_{\lambda_n}(f) \geq-\frac{c_2H(Q^\prime)}{2\sqrt{n}},
\end{eqnarray*}
where
\[
r_{n}=c_3\frac{H(Q^\prime)}{\sqrt{n}}\delta^2\|T-f^*\|^2.
\]
\end{Theorem}

\begin{pf} By Theorem
\ref{cor:OtherMinimizerInCorona}, for any $r,\lambda>0$, if $f\in F$
is such that $\mathbb{E}\mathcal{L}_\lambda(f) < r$, then
\[
\frac{r}{\lambda} > c \frac{\rho}{D}\|f-f^*\|^2,
\]
where $D$ and $\rho$ were defined in Theorem~\ref{cor:OtherMinimizerInCorona}.
Thus,
\[
\{ f\in F \dvt \mathbb{E}\mathcal{L}_\lambda(f) < r \} \subset\bigl\{ f\in
F\dvt \|f-f^*\| < c_4
\sqrt{r/\lambda} \bigr\},
\]
where $c_4=c_4(\rho,D)$. Hence, by Lemma~\ref{lemma:sym}, for $n \geq n_1$,
\begin{eqnarray*}
\mathbb{E}\sup_{\{f\in F: \mathbb{E}\mathcal{L}_\lambda(f) < r\}}
-P_n \mathcal{L}_\lambda(f)
&\leq&  c_5 \mathbb{E}\sup_{\{ f\in F: \|f-f^*\| \leq c_4 \sqrt
{r/\lambda} \}}
\Biggl|\frac{1}{n}\sum_{i=1}^n g_i (f-f^*)(X_i)\Biggr|
\\
&\leq& \frac{c_5}{\sqrt{n}} \operatorname{osc}_n\bigl(F,f^*,c_4 \sqrt{r/\lambda}\bigr)
\leq\frac{c_5}{\sqrt{n}} \eta H(Q^\prime) ,
\end{eqnarray*}
provided that $c_4 \sqrt{r/\lambda} \leq\delta$. Thus, for an
appropriate choice of $\eta$ (e.g., $\eta=c_1c_2/(4c_5)$ would do)
and setting $r_{n}:=(c_3/(2c_4^2)) n^{-1/2}H(Q^\prime)\delta^2$
(which is smaller than $\delta^2\lambda_n/c_4^2$), it is evident that
\[
\mathbb{E}\sup_{\{f\in F: \mathbb{E}\mathcal{L}_{\lambda_n}(f) <
r_{n}\}} -P_n
\mathcal{L}_{\lambda_n}(f) \leq\frac{c_1c_2}{4
\sqrt{n}}H(Q^\prime).
\]
Therefore, with $\mu^n$-probability at least $1-c_1/2$,
\[
\sup_{\{f\in F: \mathbb{E}\mathcal{L}_{\lambda_n}(f) < r_{n}\}} -P_n
\mathcal{L}_{\lambda_n}(f) \leq\frac{c_2H(Q^\prime)}{2\sqrt{n}},
\]
as claimed.
\end{pf}

We can now prove our main result.

\begin{pf*}{Proof of Theorem \ref{thm:main}} By
Theorem \ref{thm:gauss-proc} applied to the set $Q^\prime$, there
exists some integer $n_0=n_0(Q^\prime)$ such that for every $n \geq n_0$,
with $\mu^n$-probability at least $c_1$,
%
\begin{equation}\label{eq:Final1}
\inf_{\mathcal{L}(f) \in Q^\prime} P_n\mathcal{L}_{\lambda_n}(f)
\leq
-c_2\frac{H(Q^\prime)}{\sqrt{n}},
\end{equation}
where $c_1$ and $c_2$ are two absolute constants.

By Theorem~\ref{thm:UCLT}, for any integer $n\geq n_1$, with
$\mu^n$-probability at least $1-c_1/2$,
%
\begin{equation}\label{eq:Final2}
\inf_{\{f\in F: \mathbb{E}\mathcal{L}_{\lambda_n}(f) < r_{n}\}} P_n
\mathcal{L}_{\lambda_n}(f) \geq-\frac{c_2H(Q^\prime)}{2\sqrt{n}}.
\end{equation}

Hence, combining equations (\ref{eq:Final1}) and (\ref{eq:Final2}),
with $\mu^n$-probability at least $c_1/2$, the excess risk of $\hat
{f}_{\lambda_n}$ is such that $\mathbb{E}[\mathcal{L}_{\lambda
_n}(\hat{f}_{\lambda
_n})|D] \leq
-c_2H(Q^\prime)/(\sqrt{n})$, while for every function $f\in F$ with
$\mathbb{E}\mathcal{L}_{\lambda_n}(f) < r_{n}$, the empirical excess risk
satisfies $P_n \mathcal{L}_{\lambda_n}(f) \geq-c_2H(Q^\prime
)/(2\sqrt{n})$.
Therefore, the empirical risk minimization algorithm has an excess
risk (conditionally on the data $D$) larger than $r_{n}$, with
probability greater than $c_1/2$, as claimed.
\end{pf*}

\section*{Acknowledgements} This research was supported in part by
Australian Research Council Discovery Grant
DP0559465 and by Israel Science Foundation Grant 666/06.


\printhistory

\end{document}